\documentclass{amsart}
\usepackage{amsmath,amssymb,amscd}
\usepackage[all]{xy}
\newtheorem{theorem}{Theorem}[section]
\newtheorem{remark}[theorem]{Remark}
\newtheorem{proposition}[theorem]{Proposition}
\newtheorem{lemma}[theorem]{Lemma}
\newtheorem{definition}[theorem]{Definition}
\newtheorem{example}[theorem]{Example}

\newtheorem{corollary}[theorem]{Corollary}

\usepackage{pdfpages}
\usepackage{pgf,tikz}
\usetikzlibrary{arrows}
\usepackage{pdfpages}
\usepackage{pgf,tikz}
\usetikzlibrary{arrows}
\usepackage{tikz-cd}
\usepackage{graphicx}%
\usepackage{multirow}%
\usepackage{amsmath,amssymb,amsfonts}%
\usepackage{amsthm}%
\usepackage{mathrsfs}%
\usepackage[title]{appendix}%
\usepackage{xcolor}%
\usepackage{textcomp}%
\usepackage{manyfoot}%
\usepackage{booktabs}%
\usepackage{algorithm}%
\usepackage{algorithmicx}%
\usepackage{algpseudocode}%
\usepackage{listings}%

\title{Filters and congruences in weakly complemented lattices} 
\author[ Y.L. Tenkeu Jeufack]{Yannick Léa Tenkeu Jeufack}
\author[L. Kwuida]{Leonard Kwuida}
\address{University of Yaounde 1, Faculty of Sciences,
	Department of Mathematics, Laboratory of Algebra, P.O. Box 812,
	Yaounde, Republic of Cameroon} 
\email{yannick.tenkeu@facsciences-uy1.cm} 
    
\address{Bern University of Applied Sciences, Brückenstrasse 73, 3005 Bern, Switzerland.}
\email{leonard.kwuida@bfh.ch} 
\dedicatory{To the memory of Professor Marcel Tonga, 11.04.1951-15.09.2025
}
\thanks{This work was partially supported by the IMU-Simons Research Fellowship program, grant No.~640. We acknowledge the support of Bern University of Applied Sciences (BFH) during this visit.
}

\begin{document}
	\begin{abstract} 	In this paper, we show that given a weakly dicomplemented lattice (WDL)  $\mathcal{L}=(L; \vee, \wedge, ^{\Delta}, ^{\nabla}, 0, 1)$, $^{\Delta}$ induces a structure of a dual weakly complemented lattice in the lattice $(F(L), \subseteq)$ of filters of $\mathcal{L}$. We prove that the set of dense elements of $F(L)$ forms a nearlattice, and the set of principal filters of $\mathcal{L}$ forms a dual weakly complemented lattice that is  dually isomorphic to the weakly complemented lattice (WCL) $(L,\wedge, \vee, ^{\Delta}, 0, 1)$.\par  Each filter of the dual skeleton $\overline{S}(L)$ of $L$ constitutes a base of some filter in $L$,  called an S-filter, and it is proved that S-filters form a complete lattice isomorphic to the complete lattice of filters of $\overline{S}(L)$.  S-primary filters are introduced and investigated, and it is shown that there exists a bijection between the set of prime filters of $\overline{S}(L)$  and the set of S-primary filters of $\mathcal{L}$.  Furthermore, each maximal filter of a WDL $\mathcal{L}$ is a primary filter, though  there  exist primary filters of $\mathcal{L}$ that are not maximal.The congruences generated  by  filters in a distributive weakly complemented lattice are characterized. Finally,  simple and subdirectly irreducible distributive  weakly complemented lattices are also characterized using S-filters.
	\end{abstract}

	\maketitle

	\section{Introduction}
	
	Weakly dicomplemented lattices (WDLs) were introduced by Rudolf Wille  to develop a Boolean concept Logic based on concept as unit of thought \cite{RWille2}.
	Building on the formalization of the notion of concept introduced in  \cite{RWille82} and the resulting complete lattices, the meet and the join operators are used to encode the disjunction and conjunction. To introduce a negation on  concepts,  Wille follows the idea of Georges Boole in \cite{Boole}, but came out with several formalizations  that give rise to  concept algebras (with WDLs as underline algebraic structures) \cite{kwdd,RWille2} on one hand, and to protoconcepts (with double Boolean algebras as underlying structures) \cite{tenk2, RWille82} on the other. These algebras are relatively recent and have attracted the attention  of many researchers \cite{BrecknerSac, kwdd, HP2023, tenk1, tenk2, tenk3,Vorm2}.\par 
	The investigation of weakly dicomplemented lattice  has  so far focused primarily on the representation problem \cite{kw3, kwdd}, as well as on  congruences and homomorphisms \cite{kw3, kwdd, kw5}.  The prime ideal theorem \cite{kwdd} plays a central role in these attempts at  representation.  However, from a category  point of view, it seems important to explore filters and congruences  to understand the structure of the concept  algebra associated with a given context. \par 
	To contribute to a deeper understanding of WDLs, our  goal is to  investigate algebraic notions  such as filters and ideals, which are fundamentals tools for topological representation, and to continue the study of  congruence relations in special classes.  
	The introduction of new classes of filters in WDLs and the connection between the filters of the lattice $L$ and those of its dual skeleton will  provide further insights into the structure of WDLs.  This investigation on filters and ideals constitutes a preliminary work that will support  the  study of spectral theory and the  construction of  Priestley spaces \cite{HP70} for WDLs   and distributive WCLs. 
	This paper  is organized into six sections. Section 2 presents the basic definitions and algebraic properties of  WDLs; along with some useful known results on lattices.   Section 3 addresses the algebraic properties of the set  of filters of a WDL. We show that the set of filters of a WDL $\mathcal{L}$ can be endowed  with a dual weak complementation and show that a reduct $(L; \vee, \wedge, ^{\Delta}, 0, 1)$ of $L$ is dual isomorphic to the  algebra formed by the principal filters. Moreover,  we show  that  the set  of dual  dense filters of $L$ forms a nearlattice, which is distributive if $L$ is itself distributive. \par 
	In Section 4, we show that each filter of the skeleton of $L$ is a base ( see Definition \ref{def:base}) of some filter in $L$,  called an $S$-filter (see Definition \ref{def:S-filter}).  We also  observe that the classical join of two $S$-filter is not always an $S$-filter,  this observation allows us to characterize the S-filter generated  by a subset and the join of two S-filters. We further show that, for any WDL $L$,  there exists an injective, order preserving map between the filters of  $\overline{S}(L)$ and those of $L$. The isomorphism between  the lattice of filters of $\overline{S}(L)$ and the lattice of S-filters is also established.\par 
	In Section 5, we focus on prime and primary filters. We show that prime filters are precisely the primary filters in the dual skeleton  $L$, and that every  prime filter is a primary filter, although the converse is not necessarily true.  It is proved that every maximal filter is a primary filter, but the converse does not hold.
	\par In Section 6, we show that, if $\theta$ is a congruence relation on $L$, then the class of the unit is an $S$-filter, and conversely, if $F$ is an $S$-filter of $L$, and $L$ distributive, we characterize the congruence on $(L; \vee, \wedge, ^{\Delta}, 0,1)$ generated by $F$. We then characterize regular, simple, and subdirectly irreducible WCLs. We also prove that in a  distributive WCL, the congruences generated by filters commute. 
	\section{Preliminaries}
	In this section, we recall some basic definitions and results that will be needed in the subsequent sections. For preliminaries on concept algebras, see \cite{kwdd}.
	\begin{definition}\cite{kwdd}\label{theo:concept} 
		A \textbf{Weakly Dicomplemented lattice} is an algebra $(L; \wedge,\vee,^{\Delta},^{\nabla},0,1)$ of type $(2,2,1,1,0,0)$ for which $(L;\wedge,\vee,0,1)$ is a bounded lattice satisfying the following properties.\\
		
		\begin{minipage}{7cm}\begin{itemize}
				\item[(1)]  $x^{\Delta\Delta}\leq x$,
				\item[(2)] $x\leq y\implies y^{\Delta}\leq x^{\Delta}$,
				\item[(3)] $(x\wedge y)\vee(x\wedge y^{\Delta})=x$,
			\end{itemize}
		\end{minipage}\begin{minipage}{6cm}
			\begin{itemize}	\item[(1')] $x^{\nabla\nabla}\geq x$,
				\item[(2')] $x\leq y\implies y^{\nabla}\leq x^{\nabla}$,
				\item[(3')] $(x\vee y)\wedge (x\vee y^{\nabla})=x$.
			\end{itemize}
		\end{minipage}
	\end{definition}
	The reduct $(L; \vee, \wedge, ^{\Delta}, 0, 1)$ (resp. $(L; \vee, \wedge, ^{\nabla}, 0, 1)$) of $L$ is called \textbf{ weakly complemented} (resp. \textbf{dual weakly complemented}) \textbf{lattice}. $x^{\Delta}$ is called \textbf{weak complement} of $x$, and $x^{\nabla}$ is called  the \textbf{dual weak complement} of $x$. The pair $(^{\Delta}, ^{\nabla})$ is referred to as a \textbf{weak dicomplmentation}.\par 
	An \textbf{orthocomplemented lattice } is a bounded lattice with orthocomplementation, i.e a unary operation $^{\bot}: L\to L$ such that, for all $a, b\in L$:\\(i)
	$a\leq b\implies b^{\bot}\leq a^{\bot}$, 
	(ii) $(a^{\bot})^{\bot}=a$, (iii) $a\vee a^{\bot}=1$, (iv) $a\wedge a^{\bot}=0$ \cite{IvanChaj, Sylvia-Karl}).\par 
	Before providing more properties, we give some  examples of WDLs.
	\begin{enumerate}
		\item  Each Boolean algebra  can be made into a weakly dicomplemented lattices by defining $x^{\Delta}:=x', x^{\nabla}:=x'$ (where $'$ is the Boolean complement)(\cite{kwdd}).			
		\item Each bounded lattice can be endowed with a trivial weak dicomplementation by defining\\
		$a^{\Delta}=\left\{
		\begin{array}{ll}
			1 & \hbox{if}~ a\neq 1,\\
			0 & \hbox{if}~ a=0
		\end{array}\right. $ 
		$~a^{\nabla}=\left\{
		\begin{array}{ll}
			0 & \hbox{if}~ a\neq 0\\
			1, & \hbox{if}~ a=0.
		\end{array}
		\right.
		$ (\cite{kwdd}).
		\item If $L$ is a WDL and $X\neq\emptyset$, then $L^{X}$ can be  endowed with a WDL structure, where operations are defined pointwise, and constant functions are defined by $\tilde{0}(x)=0$ and $\tilde{1}(x)=1$ for all $x\in X$.
		
		\item Consider the lattice $L_{7}$, with $L_{7}=\{0,u, v, a, b, w, 1\}$ represented by the diagram in Fig 1,  with  unary operations $^{\Delta}$ and $^{\nabla}$ defined by the  Cayley  tables  given in (T: 1). The algebra $(L_{7};  \vee, \wedge, ^{\Delta}, ^{\nabla}, 0,1)$ is a weakly di-complemented lattice.\end{enumerate}
	\begin{minipage}{6cm}
		(T: 1)$\quad$\begin{tabular}{|c|c|c|c|c|c|c|c|}\hline
			$x$ & $0$ & $u$ & $v$ & $a$ & $b$ &$w$& $1$ \\
			\hline
			$x^{\Delta}$ & $1$ & $1$ & $1$ & $b$ & $a$ & $1$ &$0$ \\
			\hline
			$x^{\nabla}$ & $1$ & $v$ & $u$ & $0$ & $0$ & $0$ & $0$\\
			\hline
		\end{tabular}
	\end{minipage}
	\begin{minipage}{7cm}
		\begin{tikzpicture}[]
			\clip(-4.17,-3.7) rectangle (3.41,3);
			\draw [line width=2pt] (-2,1)-- (-1,2);
			\draw [line width=2pt] (-1,2)-- (0,1);
			\draw [line width=2pt] (0,1)-- (-1,0);
			\draw [line width=2pt] (-1,0)-- (-2,1);
			\draw [line width=2pt] (-2,1)-- (-2,-1);
			\draw [line width=2pt] (-2,-1)-- (-1,-2);
			\draw [line width=2pt] (-1,-2)-- (-1,0);
			\draw [line width=2pt] (0,1)-- (0,-1);
			\draw [line width=2pt] (0,-1)-- (-1,-2);
			\draw (-1.13,-2.16) node[anchor=north west] {0};
			\draw (-2.5,-0.68) node[anchor=north west] {$u$};
			\draw (0.17,-0.56) node[anchor=north west] {$v$};
			\draw (-1.3,0.64) node[anchor=north west] {$w$};
			\draw (-1.07,2.5) node[anchor=north west] {$1$};
			\draw (-2.47,1.44) node[anchor=north west] {$a$};
			\draw (0.1,1.26) node[anchor=north west] {$b$};
			\draw (-1.49,-2.62) node[anchor=north west] {$Fig ~1$};
	\end{tikzpicture}\end{minipage}~\\
	
	Let $(L;  \leq)$ be a partially ordered set. A \textbf{closure operator} $c$ in $L$ is a map $$c: (L,\leq )\to (L, \leq)$$
	such that for all $x, y\in L$, the following three conditions hold:
	\begin{itemize}\item[(i)] \textbf{extensive}: 
		$x\leq c(x)$,
		\item[(ii)] \textbf{isotone}:  $x\leq y\Rightarrow c(x)\leq c(y)$,
		\item[(iii)] \textbf{idempotent}:  
		$c(c(x))=c(x)$.  
	\end{itemize}
	If we replace (i) by (i)' $c(x)\leq x$ (\textbf{coextensive}), then $c$ is  called an \textbf{ interior operator} on $\textbf{L}$ (see \cite{kwdd}). In the next,  $L$ denotes a WDL $(L; \vee, \wedge, ^{\Delta}, ^{\nabla}, 0,1)$.  
	The map: $x\mapsto x^{\nabla\nabla}$ is a closure operator on $\textbf{L}$,  and 
	the map $x\mapsto x^{\Delta\Delta}$ is an interior operator on $\textbf{L}$.
	The set of closed elements $$S(\textbf{L})=\{x\in L\mid x^{\nabla\nabla}=x\}$$ is called \textbf{skeleton} of $\textbf{L}$,  and the set of interior elements  $$\overline{S}(\textbf{L})=\{x\in L\mid x^{\Delta\Delta}=x\}$$  is called the \textbf{dual skeleton}.\par 
	
	We define the operations $\sqcup$ and $\overline{\sqcap}$ on $L$ by
	$$  x\sqcup y:=(x^{\nabla}\wedge y^{\nabla})^{\nabla}~\text{,}~ x\overline{\sqcap} y:=(x^{\Delta}\vee y^{\Delta})^{\Delta},\quad x\underline{\sqcup} y=(x\vee y)^{\Delta\Delta}.$$
	According to 
	(\cite{kwdd}, Prop. 1.3.4)\label{prop:ortho}, the algebras 
	$$(S(\textbf{L}); \wedge, \sqcup,^{\nabla},  0, 1)~\text{ and}~ (\overline{S}(\textbf{L}); \overline{\sqcap}, \vee,^{\Delta}, 0, 1)$$ are ortholattices. A nonempty subset  $F$ of   $\textbf{L}$ is called \textbf{ filter} if (i) $x\leq y$ and $x\in F$ implies $y\in F$ and  (ii)  
	$x,y\in F\Rightarrow x\wedge y\in F, \forall x,y\in L.$ 
	Dually, one defines an ideal. A proper filter (resp. ideal) is a filter (resp. ideal) that is different from $L$.
	A filter $F$ is called  a \textbf{primary filter} (rep. \textbf{prime filter}) if: for all $x\in L$, $x\in F$ or  $x^{\Delta}\in F$ (resp. $x\vee y\in F \Rightarrow x\in F$ or $y\in F$, $\forall x,y\in L$).
	The filter $[a)=\{x\in L\mid a\leq x\}$ is called the principal  filter generated by $a\in L$. Let  $F(\textbf{L})$ be the set of  filters of $\textbf{L}$. Ordered under set  inclusion,  $(F(L), \subseteq )$ is a complete lattice, with  least element $\{1\}$ and greatest element $L$.\par  Clearly, the meet of two filters $F$ and $G$ is given by set intersection and the join  is denoted by $F\vee G$.
	From now on, let $L$ denote a nontrivial  WDL  and let $L^{\Delta}$ denote its reduct  $(L; \vee, \wedge, ^{\Delta}, 0,1)$ which is a WCL. The following identities will be  useful in calculations:
	\begin{proposition} \cite{kwdd}.\label{prop: axiom} Let $x, y\in L$. The following statement holds. \vspace{0.3cm}
		
		\begin{minipage}{6cm}
			\begin{itemize}
				\item[(4)] $(x\wedge y)^{\Delta}=x^{\Delta}\vee y^{\Delta}$.
				\item[(5)] $x\overline{\sqcap} y=(x\wedge y)^{\Delta\Delta}$ 
				\item[(6)] $x^{\Delta}\leq y\Leftrightarrow y^{\Delta}\leq x$

				\item[(7)] $x^{\Delta}\leq y\Rightarrow y\vee x=1$.
				\item[(8)] $y\wedge x=0\Rightarrow y^{\Delta}\geq x$.
				
				\item[(9)] 
				$x\vee x^{\Delta}=1$.
				
			\end{itemize}
		\end{minipage}
		\begin{minipage}{7cm}
			\begin{itemize}
				\item[(4')] $(x\vee y)^{\nabla}=x^{\nabla}\wedge y^{\nabla}$.
				\item[(5')]  $x\sqcup y=(x\vee y)^{\nabla\nabla}$.
				\item[(6')] $x^{\nabla}\geq y\Leftrightarrow y^{\nabla}\geq x$.
				\item[(7')] $x^{\nabla}\geq y\Rightarrow y\wedge x=0$.
				\item[(8')] $y\vee x=1\Rightarrow y^{\nabla}\leq x$.
				\item[(9')] a) $x\wedge x^{\nabla}=0$. b) $x^{\nabla}\leq x^{\Delta}$.
			\end{itemize}
		\end{minipage}
	\end{proposition}
	\begin{remark}\label{lem:comp}
		Let $x, y\in L$, the following statements hold:
		\begin{enumerate}
			\item (i) $x\overline{\sqcap} y\leq x\wedge y$, (ii)   $x\underline{\sqcup} y\leq x\vee y\leq x\sqcup y$.
			\item (i)  $a\mapsto x\sqcup a$ preserves $\leq$ and $\sqcup$. (ii) $x\mapsto x\overline{\sqcap} a$ preserves $\leq$ and $\overline{\sqcap}$. 
			\item (i) $x\vee(x\overline{\sqcap} y)=x$, (ii) $x\overline{\sqcap} (x\vee y)=x^{\Delta\Delta}$.
			\item (i) $x\wedge (x\sqcup y)=x$, (ii) $x\sqcup (x\wedge y)=x^{\nabla\nabla}$.  
			\item (i) $x^{\Delta\Delta}\overline{\sqcap} y^{\Delta\Delta}=x^{\Delta\Delta}\overline{\sqcap} y=x\overline{\sqcap} y$. (ii) $x^{\nabla\nabla}\sqcup y^{\nabla\nabla}=x^{\nabla\nabla}\sqcup y=x\sqcup y$.
			\item If $x\leq y$ and $a\leq b$, then (i) $x\overline{\sqcap} a\leq y\overline{\sqcap} b$, (ii) $x\sqcup a\leq y\sqcup b$.
		\end{enumerate}
	\end{remark}
	We recall the well-known characterization of the join of two filters in the following proposition. Note that for any nonempty subset $X$ of $L$, $[X)$ denotes the filter of $L$  generated by $X$.
	\begin{proposition} \label{prop: engendrer} For any  $G, F\in F(\textbf{L})$,  $a, b\in L$ and $X\subseteq L$, the following statements hold:
		\begin{enumerate}
			\item   $G\vee F=\{x\in L\mid g\wedge f\leq x, \text{for some}~g\in G, f\in F\}$.
			\item $[X)=\{x\in L\mid \exists u_{1}, \ldots, u_{n}\in X, u_{1}\wedge \ldots\wedge u_{n}\leq x, n\geq 1\}.$
			\item  (i) $[a)\vee[b)=[a\wedge b)$, (ii) $[a)\cap[b)=[a\vee b)$.
		\end{enumerate}
	\end{proposition}	
	Note that in Boolean algebras, filters and ideals are dual notions. However, in a WDL, filters and ideals do not form such pairs; nevertheless, each result concerning filters has a corresponding dual result for ideals. We introduce certain algebraic structures to provide a more detailed description of the set of filters (resp. ideals) of a WDL. We also investigate specific elements in WDLs, namely dense and codense elements. To this end, we introduce nearlattices and dual nearlattices.
	
	\begin{definition}(\cite{IvanChaj})\begin{enumerate}
			\item A \textbf{join semi-lattice with greatest element} is an algebra $(A;\vee, 1)$ of type $(2,0)$ such that the operation $\vee$ is idempotent, commutative, associative and $x\vee 1=1$ for all $x\in A$. 
			\item A \textbf{meet semilattice} is defined dually.
			\item A \textbf{nearlattice} is a join  semilattice $A$ such that for each $a\in A$,  the principal filter $[a)$ is a bounded lattice.
		\end{enumerate} 
	\end{definition}

	\begin{definition} An element $x$ of $L$ is said to be \textbf{dense} (resp. \textbf{codense}),  if $x^{\nabla}=0$ (resp. $x^{\Delta}=1$).
	\end{definition}
	
	Consider the set $D(L)=\{x\in L\mid x^{\nabla}=0\}$ and $ \overline{D}(L)=\{x\in L\mid x^{\Delta}=1\}$.
	\begin{remark}\label{prop: semilattices}  The following statements hold:\vspace{0.2cm}\\
		\begin{minipage}{7.5cm}
			\begin{itemize}
				\item[(i)]  $1\in D(L)$.
				\item[(ii)]  $D(L)$ is an order filter of $L$.
				\item[(iii)]  If $x, y\in D(L)$, then $x\vee y\in D(L)$.
			\end{itemize}
		\end{minipage}
		\begin{minipage}{7.5cm}
			\begin{itemize}
				\item[(i')] $0\in \overline{D}(L)$.
				\item[(ii')]$ \overline{D}(L)$ is an order ideal of $L$. 
				\item[(iii')]  If $x, y\in \overline{D}(L)$, then $x\wedge y\in \overline{D}(L)$. 
			\end{itemize}
		\end{minipage}
	\end{remark}	
	\begin{proposition}\label{theo:near}
		The following statements hold:  
		$D(L)$ (resp. $\overline{D}(L)$) is a nearlattice (resp. a dual nearlattice). Moreover, it is distributive if $L$ is distributive.
	\end{proposition}
	
	\begin{proof}
		Let $a \in D(L)$,  $[a)$ is a bounded lattice. Let $x, y \in D(L)$. Then, by (ii) of Remark  \ref{prop: semilattices}, one has $x \vee y \in D(L)$. Since $\vee$ is commutative, associative, and idempotent on $L$, and $1 \in D(L)$, it follows that $D(L)$ is a join semilattice. Therefore, $(D(L), \vee, 1)$ is a nearlattice.  
		A similar argument applies to $\overline{D}(L)$.
	\end{proof}
	
	\begin{remark}
		If $D(\mathbf{L})$ has a least element, then $D(\mathbf{L})$ is a filter.
	\end{remark}
	
	\section{The structure of the set $F(L)$} 
	
	In this section, we show that the set $F(L)$ carries the  structure of a dual weakly complemented lattice, and  that its dense elements form a nearlattice.\par
	
	For  $ X, Y\subseteq L$,  define  the sets $X^{\star}$ and $X^{+}$ as follows:  \begin{itemize}
		\item[(1)] (i)  $X^{\star}:=\{a\in L\mid \forall x\in X, x^{\Delta}\leq a\}$. (ii)  $X^{+}:=\{x\in L\mid \forall a\in X, x\vee a=1\}$.
	\end{itemize}
	
	Moreover,  for any $X\subseteq \overline{S}(L)$, set 
	$X^{\overline{\star}}:=\{x\in \overline{S}(L)\mid \forall a\in X, a^{\Delta}\leq x\}.$
	\begin{lemma}\label{lem:sferme} The following statements hold:

		\begin{enumerate}
			\item (i) If $X\subseteq Y$, then $Y^{\star}\subseteq X^{\star}$.
			(ii) $X\subseteq X^{\star\star}$.  (iii)  (a) $[ X)^{\star}=X^{\star}$.
			(iv)  $X^{\star\star\star}=X^{\star}$. 
			\item  $X^{\star}$ is a filter of $L$. Moreover,   $L^{\star}=\{1\}$ and $ \{1\}^{\star}=L$.
			\item The map $X\mapsto X^{\star\star}$ is a closure operation on $F(L)$.
			\item The map  $X\mapsto X^{\overline{\star}\overline{\star}}$ is a closure operation on $F(\overline{S}(L))$.  
		\end{enumerate}
		
	\end{lemma}
	
	\begin{proof}
		For item (1), (i) is straightforware to check. For (ii),  let $x\in X$ and  $y\in X^{\star}$. By definition of $^{\star}$, we have $x^{\Delta}\leq y$. Moreover, by  (1) of Definition \ref{theo:concept}, 
		$y^{\Delta}\leq x^{\Delta\Delta}\leq x,$  so $x\in X^{\star\star}$. Hence  $X\subseteq X^{\star\star}$.
		Now let us prove (iii).  Since $X\subseteq [X)$, it follows from (1)(i) that  $[X)^{\star}\subseteq X^{\star}$. For the reverse inclusion,  let $a\in X^{\star}$. To show that $a\in [X)^{\star}$, let $u\in [X)$. By (2) of Proposition \ref{prop: engendrer}, there exists $x_{1}, \ldots, x_{n}\in X$ such that $$x_{1}\overline{\sqcap}\ldots\overline{\sqcap} x_{n}\leq x_{1}\wedge \ldots\wedge x_{n}\leq u.$$ 
		
		Since  $x_{i}^{\Delta}\leq a$ for all $1\leq i\leq n$, we obtain $$u^{\Delta}\leq (x_{1}\overline{\sqcap}\ldots\overline{\sqcap} x_{n})^{\Delta}=x_{1}^{\Delta}\vee\ldots\vee x_{n}^{\Delta}\leq a,$$ 
		which shows that  $a\in [X)^{\star}$. Therefore,  $X^{\star}\subseteq [X)^{\star}$ and hence $X^{\star}=[X)^{\star}$.\\
		For (iv), it is clear that  $X^{\star}\subseteq X^{\star\star\star}$ by (ii). Applying (i) to (ii), we obtain  $X^{\star\star\star}\subseteq X^{\star}$. Thus $X^{\star\star\star}=X^{\star}$ and (iv) follows.\par 
		For item (2):\\
		(i) Assume $X\neq \emptyset$. Then $x^{\Delta}\leq 1$ for all $x\in X$, so $1\in X^{\star}$. Let $a,b \in X^{\star}$ and $x\in X$. By  assumption,  $x^{\Delta}\leq a, b$. In addition,  $x^{\Delta}\overline{\sqcap} x^{\Delta}=x^{\Delta}\leq a\overline{\sqcap} b\leq a\wedge b,$ so $a\overline{\sqcap} b, a\wedge b\in X^{\star}$. Now  assume  $a\leq b$ and $a\in X^{\star}$. For  $x\in X$,  we have  $x^{\Delta}\leq a\leq b$, so  $b\in X^{\star}$. Therefore  $X^{\star}$ is a filter of $L$ and   it is closed under  $\overline{\sqcap}$. \\(ii) is obvious.
		The proof of (3) follows from (1),  and the proof of (4)  is similar to that of (2).
		
	\end{proof}
	Clearly $X^{+}$ is a filter of $L$ when $L$ is distributive.
	From Proposition \ref{lem:sferme} we can work with filters of $L$ (resp. $\overline{S}(L)$) instead of arbitrary subset of $L$  when using the operations $^{\star}$ and $^{\overline{\star}}$. 
	\begin{proposition}\label{prop: filpseudo}  The map $^{\star}$ is a dual weak complementation in $F(L)$.		
	\end{proposition}
	\begin{proof}
		(1) Let $F$ be a filter of $L$ and $F^{\star}$ as above. Then  $F^{\star}$ is a filter of $L$ by (2)(i) of Lemma  \ref{lem:sferme}.   It remains to show that $^{\star}$ verifies  conditions (1'), (2') and (3') in Definition \ref{theo:concept}. Let $F$ and $G$ be two filters of $L$.\\	
		(i) By (i) of  Lemma \ref{lem:sferme}, one has: $F\subseteq G$ implies  $G^{\star}\subseteq F^{\star}$. Therefore,  (2')  holds.\\	
		(ii) For (3'), we show that $(F\vee G)\cap (F\vee G^{\star})=F$. Clearly $F\subseteq (F\vee G)\cap (F\vee G^{\star})$.  For the reverse inclusion, let $z\in (F\vee G)\cap (F\vee G^{\star})$. 
		Then $z\in F\vee G$ and $z\in F\vee G^{\star}$ and by  (1) of Proposition \ref{prop: engendrer},  there  exist $u, a\in F$, $v\in G$ and $b\in G^{\star}$ such that $$u\wedge v\leq z~(1), \quad \quad ~\quad a\wedge b\leq z~(2).$$
		Since  $b\in G^{\star}$ and $v\in G$, one has $v^{\Delta}\leq b$ (3).
		From  (3), we obtain  $(a\wedge u)\wedge v^{\Delta}\leq (a\wedge u)\wedge b\leq z~ (4),$ and in addidition, $(a\wedge u)\wedge v\leq u\wedge v\leq z \quad (5).$ 
		From (4) and (5), we get   $$[(a\wedge u)\wedge v]\vee[(a\wedge u)\wedge v^{\Delta}]\leq z.$$ 
		Using (3) of Definition \ref{theo:concept}, we deduce that   $a\wedge u\leq z$. Since   $a\wedge u\in F$,  it follows that $z\in F$. 
		Therefore,  the desired equality holds.\\
		(iii) It remains to show that $F\subseteq F^{\star\star}$. This holds by (1) (ii) of Lemma  \ref{lem:sferme}. Thus (1') of Definition \ref{theo:concept} is satisfied,  and the proof  is complete.
	\end{proof}
	Now we can  state our  main result of this  section.
	\begin{theorem}\label{theo:dwdl}   The algebra  $(F(L); \cap, \vee, ^{\star}, \{1\}, L)$ is a dual weakly complemented lattice.
	\end{theorem}
	\begin{proof}
		(1) It is clear that $(F(L), \cap, \vee, \{1\}, L)$ is a complete lattice and hence  a bounded lattice. By  Proposition \ref{prop: filpseudo}, the operation $^{\star}$ is a dual weak complementation on $F(L)$. Therefore,  the algebra $(F(L); \cap, \vee,^{\star}, \{1\}, L)$ is a dual weakly complemented lattice. 
	\end{proof}

	Consider the condition: $(\star)~ x\vee y=1\Rightarrow x^{\Delta}\leq y$.
	\begin{corollary} The following statements hold:
		\begin{itemize}\item[(i)] If $(L; \wedge , \vee, ^{\Delta}, 0, 1)$ is distributive and verifies the condition ($\star$), then $(F(L); \cap, \vee, ^{\star}, \{1\}, L)$ is a pseudocomplemented lattice verifying $F^{\star}=F^{+}$ for all $F\in F(L)$. Moreover it is distributive if $L$ is distributive.  
			\item[(ii)] 
			For all $F, G\in F(L)$,  if $F\subseteq G$,  then $G\cap (G\vee F^{\star})=G$.\end{itemize}
	\end{corollary}
	
	\begin{proof}
		(i) Assume that $L$ is distributive. Then $F(L)$ is clearly a distributive lattice.  By  Theorem \ref{theo:dwdl},  $(F(L); \cap, \vee, ^{\star}, \{1\}, L)$ is a dual weakly complemented lattice. For any filter $F$ of $L, F\cap F^{\star}=\{1\}$.  Assume that condition $(\star$) holds for $(L; \wedge, \vee, ^{\Delta}, 0,1)$. Let $F, G$ be two filters of $L$ such that $G\cap F=\{1\}$. We show that $G\subseteq F^{\star}$.\\
		Let $z\in G$ and $u\in F$. Then $z\vee u\in F\cap G=\{1\}$, hence   $z\vee u=1$. By assumption, this implies   $u^{\Delta}\leq z$, so $z\in F^{\star}$. Therefore $G\subseteq F^{\star}$, and the claim follows. \par 
		It remains to show that $F^{+}=F^{\star}$. Using (9) of Proposition \ref{prop: axiom}, we obtain  $F^{\star}\subseteq F^{+}$. From the  definition of $F^{\star}$ and $F^{+}$,  combined with the assumption,  we obtain $F^{+}\subseteq F^{\star}$. Hence  $F^{+}=F^{\star}$ whenever $L$ satisfied the implication $(\star$).\\
		(ii) This is straighforward.
	\end{proof} 
	
	\begin{proposition} If $L$ is a distributive lattice, then:
		$(F(L); \wedge , \vee, ^{+}, \{1\}, L)$ is a pseudocomplemented distributive  lattice.
	\end{proposition}
	\begin{proof}
		Assume that $L$ is a distributive WDL. It is easy to show that $F(L)$ is a distributive bounded lattice  and for any $F\in F(L)$, $F^{+}$ is  a pseudocomplement of $F$ in $F(L)$. 
	\end{proof}
	Let $D(F(L)):=\{F\in F(L)\mid F^{\star}=\{1\}\}$.

	\begin{corollary}
		\begin{enumerate}
			\item The algebra $(S(F(L)); \wedge, \sqcup, {^\star}, \{1\}, L)$ is a complete ortholattice, where the meet is given by 
			$F \wedge G = F \cap G, \quad \text{and the join by } F \sqcup G = (F \vee G)^{\star\star}.$
			\item The algebra $(D(F(L)); \sqcup, L)$ is a nearlattice; moreover, it is distributive whenever $L$ is distributive.
		\end{enumerate}
	\end{corollary}
	
	\begin{proof}
		(1) Since $^{\star}$ is a dual weak complementation on $F(L)$, and $^{\star\star}$ is a closure operation, it follows that $(S_{F}(L); \wedge, \sqcup)$ is a complete lattice with the join and meet defined as stated. Hence, it is an ortholattice (see \cite{kwdd}, Proposition \ref{prop:ortho}). \\[6pt]
		(2) This follows from Theorem \ref{theo:dwdl},   Proposition \ref{theo:near}, and the fact that $F(L)$ is a dual weakly complemented lattice.
	\end{proof}
	
	Next, we consider the particular case of principal filters. Let $F_{p}(L)$  denotes the set of principal filters of $L$. It $L$ is a Boolean algebra,  then $F_{p}(L)$ forms a Boolean  algebra isomorphic to $L$.  In the following, we show that a similar result holds for distributive WCL. \par 
	For  $a\in L$, we have   $[a)^{\star}=\{x\in L\mid a^{\Delta}\leq x\},$ is called the \textbf{dual weak annihilator} of $a$. \begin{remark}\label{prop:pclosed} Let $F\in F(L)$, $a, b\in L$. The following statements hold.
		\begin{enumerate}
			\item (i)  $[a)^{\star}=[a^{\Delta})$.  (ii) $[1)^{\star}=L$ and $[0)^{\star}=\{1\}$.
			(iii) $F^{\star}=\underset{x\in F}{\cap}[x^{\Delta})$.
			\item If $a\leq b$, then $[a)^{\star}\subseteq[b)^{\star}$.
			
		\end{enumerate}
	\end{remark}

	Let  $\Lambda(L)=\{[a)^{\star}\mid a\in L\}$.
	\begin{theorem} \begin{enumerate}\item The algebra $(F_{p}(L); \wedge, \vee, ^{\star}, [0),[1))$  is a dual weakly complemented lattice dual isomorphic to $(L; \vee, \wedge, ^{\Delta}, 0, 1)$.
			\item The  algebra $(\Lambda(L); \wedge, \underline{\vee}, ^{\star}, \{1\}, L)$  is an ortholattice isomorphic to $(\overline{S}(L); \overline{\sqcap}, \vee, ^{\Delta}, 0, 1)$.\end{enumerate}
	\end{theorem}

	\begin{proof}	(1) \\
		(i)  By (3) of Proposition \ref{prop: engendrer}, $F_{p}(L)$ is stable with $\vee$ and $\wedge$, and by (1) of Proposition \ref{prop:pclosed} $F_{p}$ is stable with $^{\star}$. Therefore, $(F_{p}(L), \wedge, \vee,^{\star}, [1),[0))$ is a subalgebra of $(F(L); \wedge, \vee,^{\star}, \{1\}, L)$. Hence it is a dual weakly complemented lattice.\\ 
		(ii) Define 
		$\eta: L\to F_{p}(L), a\mapsto \eta(a)=[a)$.
		Clearly, $\eta$ is well defined. We have      $$\eta(a\vee b)=[a\vee b)=[a)\cap[b)=\eta(a)\wedge \eta(b)\quad \text{ and }~  \eta(a\wedge b)=[a\wedge b)=[a)\vee[b)=\eta(a)\vee\eta(b)$$ (by  (2) of Proposition \ref{prop: engendrer}). Moreover,  by (1) of Proposition \ref{prop:pclosed}  $\eta(a^{\Delta})=\eta(a)^{\star}$. Therefore,  $\eta$ is a dual homomorphism between the  weakly complemented lattice $(L; \wedge, \vee, ^{\Delta}, 0, 1)$ and the dual-weakly complemented lattice $(F_{p}; \wedge, \vee, ^{\star},\{1\}, L)$.\par
		(iii)  By definition,  $\eta$ is surjective. For injectivity, 
		assume that $\eta(a)=\eta(b)$. Then $[a)=[b)$, so $a\leq b\leq a$. By anti-symetry of $\leq $, we get  $a=b$. Hence  $\eta$ is a dual isomorphism. \par
		(2) Since $\eta$ is an isomorphism and $\overline{S}(L)$ is a dual skeleton of $L^{\Delta}$, we have  $\eta(\overline{S}(L))=\Lambda(L)$, which forms an ortholattice.   As  $\eta$ is injective, the  isomorphism is obtained by restricting  $\eta$ to $\overline{S}(L)$. 
	\end{proof}	
	\section{Coorespondence between $F(\overline{S}(L))$ and $F(L)$.} 
	
	In this section,  we investigate the  relationship between $F(\overline{S}(\textbf{L}))$  and $F(\textbf{L})$. We first introduce  a new classe of filters in $L$  and derive some of their properties.\par 
	For a given filter $F$ on $L$, we consider the following three  conditions: $$\dagger:x, y\in F\Rightarrow x\overline{\sqcap} y\in F.$$
	$$(\ddagger)~\quad \exists G\in  F(\overline{S}(L)),  F=\{x\in L\mid x^{\Delta\Delta}\in G\}.$$
	$$(\dagger\ddagger)~\quad \exists G\in  F(\overline{S}(L),  F=\{x\in L\mid \exists u\in G, u\leq x\}.$$ 
	We now present the  following useful example of WDL.
	
	\begin{example}
		Let $L_{6}:=\{0, u, v, a, b, 1\}$ be a bounded latice  with the Hasse diagram labelled in Fig 2.
		Consider the unary operations $^{\Delta}$ and $^{\nabla}$ defined by the  Cayley tables given below. One can verify that  $(L_{6}; \vee, \wedge, ^{\Delta}, ^{\nabla}, 0,1)$ is a distributive WDL,  and its reduct $L^{\Delta}$ is a distributive WCL.
	\end{example} 
	\begin{minipage}{7cm}
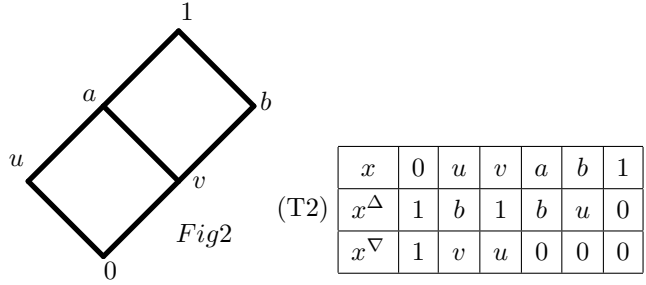
		\begin{tikzpicture}[line cap=round,line join=round,>=triangle 45,x=1cm,y=1cm]
			\clip(0.3,-2.16) rectangle (7.3, 3.4);
			\draw [line width=2pt] (4,1)-- (5,2);
			\draw [line width=2pt] (5,2)-- (6,3);
			\draw [line width=2pt] (6,3)-- (7,2);
			\draw [line width=2pt] (7,2)-- (6,1);
			\draw [line width=2pt] (6,1)-- (5,2);
			\draw [line width=2pt] (6,1)-- (5,0);
			\draw [line width=2pt] (5,0)-- (4,1);
			\draw (5.9,3.5) node[anchor=north west] {$1$};
			\draw (6.96,2.32) node[anchor=north west] {$b$};
			\draw (4.6,2.34) node[anchor=north west] {$a$};
			\draw (3.62,1.48) node[anchor=north west] {$u$};
			\draw (6.06,1.18) node[anchor=north west] {$v$};
			\draw (4.88,0.05) node[anchor=north west] {$0$};
			\draw (5.84,0.56) node[anchor=north west] {$Fig2$};
		\end{tikzpicture}
	\end{minipage}\begin{minipage}{7cm}(T2)~\begin{tabular}{|c|c|c|c|c|c|c|}
			\hline
			$x$ & $0$ & $u$ & $v$ & $a$ & $b$ & $1$ \\
			\hline
			$x^{\Delta}$ & $1$ & $b$ & $1$ & $b$ & $u$ &$0$ \\
			\hline
			$x^{\nabla}$ & $1$ & $v$ & $u$ & $0$ & $0$ & $0$\\
			\hline
		\end{tabular} 
	\end{minipage}~\\	
	One has:
	\begin{itemize}
		\item $F(\overline{S}(L_{6})) = \{\{1\}, E_{1} := \{u, 1\}, E_{2} :=\{b, 1\}, \overline{S}(L_{6})\}$.
		\item $F(L_{6}) = \{\{1\}, F_{2} := \{a, 1\}, F_{3}: = \{b, 1\}, F_{4} := \{u, a, 1\}, F_{5} := \{v, a, b, 1\}, L_{6}\}$.
		\item $F_{5}$ does not satisfy $(\ddagger)$. In fact, $a, b \in F_{5}$, but $a \overline{\sqcap} b = 0 \notin F_{5}$.
		\item Any $F \in \{F_{2}, F_{3}, F_{4}\}$ satisfies the condition $(\dagger)$.
		\item $E_{1}$ and $\overline{S}(L_{6})$ are filters of $\overline{S}(L_{6})$ that are not filters of $L_{6}$.
		\item $F_{4} = \{x \in L_{6} \mid \exists u \in E_{1}, u \le x\}$, and $F_{3} =\{x\in L\mid \exists u\in E_{2}, u\leq x\}=E_{2}$.
	\end{itemize}
	
	Note that $x\wedge x^{\Delta}$ can be different from $0$, and   $x\vee x^{\nabla}$ can be different from $1$ in $L$ (\cite{kwdd}).  For example, in    $L_{6}$, one has  $a^{\Delta}=b$ and $a\wedge a^{\Delta}=v\neq 0$. 
	
	\begin{proposition}\label{prop:equiv}
		The conditions $(\dagger)$, $(\ddagger)$, and $(\dagger\ddagger)$ are equivalent for any $F \in F(L)$.
	\end{proposition}
	\begin{proof}
		($\ddagger)\Rightarrow (\dagger\ddagger$). Assume that $F$ satisfies ($\ddagger$). Then there exists $G\in F(\overline{S}(L))$ such $F=F_{G}$. Since $G\subseteq \overline{S}(L)$, it follows immediately that $F$ also verifies  $(\dagger\ddagger$). 
		Hence $((\ddagger)\Rightarrow (\dagger\ddagger$)) holds.\\ 
		$(\dagger\ddagger)\Rightarrow (\dagger)$.  Assume that $F$ verifies $(\dagger\ddagger$). Let $x, y\in F$. By assuption, there exist $u, v\in G$ such that $u\leq x$ and $v\leq y$. By the compatibility of $\overline{\sqcap}$ with the order relation $\leq$, we have 
		$$u\overline{\sqcap} v\leq x\overline{\sqcap} y\quad \text{ and}~\quad  u\overline{\sqcap} v\in G. $$
		Thus $x\overline{\sqcap} y \in F$, which shows that  $(\dagger$) holds. Therefore, $(\dagger\ddagger) \Rightarrow (\dagger) $ holds.\\
		($\dagger)\Rightarrow (\ddagger$). Suppose that ($\dagger$) holds. We show that ($\ddagger$)  also holds. Define  $G=F\cap \overline{S}(L)$. It is easy to show that $G$ is a filter of $\overline{S}(L)$, and that $F=\{x\in L\mid x^{\Delta\Delta}\in G\}$, i. e $F$ satisfies $(\ddagger)$.\\ Therefore, the three conditions are equivalent.
	\end{proof}
	\begin{definition}\label{def:S-filter} \begin{enumerate}\item A filter $F$ of $L$ is called \textbf{S-filter} if it verifies one of the equivalent conditions $(\dagger$), ($\ddagger$) and $(\dagger\ddagger$).
			\item 
			A filter $G$ of $\overline{S}(L)$ is called \textbf{primary} if for any $x\in \overline{S}(L)$, $x^{\Delta}\in G$ or $x\in G$.
		\end{enumerate}
	\end{definition} 
	
	\begin{remark}\label{rk:join SF}
		\begin{enumerate}
			\item The filter $F_{5}=\{v, a, b, 1\}$ is not an S-filter in $L_{6}$, then a filter of $L$ is not necessarily an S-filter.
			\item  Consider $F_{2}=\{a, 1\}$ and $F_{3}=\{b, 1\}$. Both are S-filters in $L_{6}$, but their join 
			$
			F_{2}\vee F_{3}=F_{5}
			$
			is not an S-filter of $L_{6}$, where $\vee$ denotes the classical join in $F(L)$. Then, the join of two S-filters (taken as filters) is not always an S-filter.
		\end{enumerate}
	\end{remark}
	We denote by $SF(L)$ the set of S-filters of $L$.
	Now we investigate relationship between the lattice $F(\overline{S}(L))$ and $F(L)$.
	Let $F\in SF(L)$,  $G\in F(\overline{S}(L))$, set $$E:=F\cap\overline{S}(L), \quad F_{G}:=\{x\in L\mid x^{\Delta\Delta}\in G\}.$$
	\begin{proposition}\label{prop: S-filt}  The following statements hold.
		\begin{enumerate}
			\item[(i)]  $F_{G}$ is an $S$-filter of $L$ and $F_{G}\cap\overline{S}(L)=G$. \item[(ii)]$E$ is a filter of $\overline{S}(L)$  and  $F_{E}=F$.
			\item[(iii)] $F$ is a primary $S$-filter of $L$ if and only if $E$ is a primary filter of $\overline{S}(L)$.
		\end{enumerate}
	\end{proposition}
	
	\begin{proof}
		(i) We only need to show that $F_{G}$ is closed under $\overline{\sqcap}$. Let $x, y\in F_{G}$.  By the definition of $F_{G}$,  $x^{\Delta\Delta}, y^{\Delta\Delta}\in G$. Since $x^{\Delta\Delta}\overline{\sqcap} y^{\Delta\Delta}\in G$ and $x^{\Delta\Delta}\overline{\sqcap}y^{\Delta\Delta}\leq x\overline{\sqcap} y$, it follows that $x\overline{\sqcap} y\in F_{G}$. Hence $F_{G}$ is an S-filter of $L$. The second part is immediate, since  $x^{\Delta\Delta}=x$ for all $x\in G$. 
		\par 
		(ii)  It is clear that $E$ is a filter of $\overline{S}(L)$ and that $F_{E}\subseteq F$.  Let $x\in F$.  Then $x^{\Delta\Delta}\in F$ (because $F$ is an $S$-filter). Thus  $x^{\Delta\Delta}\in E\subseteq F_{E}$. Since $x^{\Delta\Delta}\leq x$, we deduce  $x\in F_{E}$. Therefore   $F\subseteq F_{E}$ and hence $F=F_{E}$.\par  	
		(iii)	Suppose  $F$ is a primary S-filter of $L$. Let $x\in \overline{S}(L)$, as $F$ is a primary filter of $L$, we have   $x\in F$ or $x^{\Delta}\in F$.  If $x^{\Delta}\in F$, then $x^{\Delta}\in E$. Otherwise, if 
		$x^{\Delta}\not\in F$, then $x\in F$, hence $x^{\Delta\Delta}=x\in F$, so $x\in E$. Thus,  $E$ is a primary filter of $\overline{S}(L)$.\par  Conversely, assume that $E$ is a primary filter of $\overline{S}(L)$. Let $x\in L$. Then $x^{\Delta}\in \overline{S}(L)$. If $x^{\Delta}\in E\subseteq F$, then $x^{\Delta}\in F$. Otherwise,   $x^{\Delta\Delta}\in E\subseteq F$, and since  $x^{\Delta\Delta}\leq x$, we deduce that $x\in F$. Therefore, $F$ is a primary filter of $L$.  
	\end{proof}
	
	\begin{definition}\label{def:base}
		A base for the filter $F$ is a subset  $F_{0}$ ($\subseteq L$) such that  $$F=\{x\in L\mid \exists u\in F_{0}, u\leq x\}.$$ 
	\end{definition}

	\begin{corollary}
		The set of $S$-filters  of $L$ forms a complete lattice.
	\end{corollary}
	
	\begin{remark}\label{prop: extfilter} Let $F\in F(L)$ and  $G_{1}, G_{2}\in F(\overline{S}(L)$. The following statements hold:
		\begin{enumerate}
			\item[(i)] $G_{1}\subseteq G_{2}\Leftrightarrow  F_{G_{1}}\subseteq F_{G_{2}}$. 
			\item[(ii)]  If  $F\cap \overline{S}(L)=G$ is a filter of $\overline{S}(L)$,  then $F_{G}\subseteq F$.   The inclusion may be strict.
			\item[(iii)] $F_{\overline{S}(L)}=L$ and $F_{\{1\}}=\{1\}$. 
			\item[(iv)] Every filter $G$ of $\overline{S}(L)$ is a base of some filter $F$ of $\textbf{L}$.
			\item[(v)]  $F \cap \overline{S}(L)$ is not necessarily a filter of $\overline{S}(L)$.   
		\end{enumerate}
	\end{remark}

	\begin{proof}
		(i) is straightforward.\\
		(ii) Suppose $F\cap\overline{S}(L)=G$ is a filter of $\overline{S}(L)$. Let $x\in F_{G}$. Then $x^{\Delta\Delta}\in G$, hence $x^{\Delta\Delta}\in F$. Since  $x^{\Delta\Delta}\leq x$, it follows that  $x\in F$. Thus $F_{G}\subseteq F$. Taking $F=\{v,a,b, 1\}$ in $L_{6}$,  one has $F\cap \overline{S}(\textbf{L}_{6})=\{b,1\}=G$ is a filter of $\overline{S}(\textbf{L})$, but $F_{G}=G\neq F$.\par
		(iii)  follows directly from the definitions. \\
		(iv) This follows from the fact that $F_{G}$ is an $S$-filter. \\
		(v) Consider the WDL $L_{7}$ with $\overline{S}(L) = \{0,1,a,b\}$, and let $F: = \{w,a,b,1\}$. Then  
		$
		\overline{S}(L_{7}) \cap F = \{a,b,1\} = G.
		$ 
		Now, $a,b \in G$, but  
		$
		a \,\overline{\sqcap}\, b = w^{\Delta\Delta} = 0 \notin G,
		$ 
		so $G$ is not a filter of $\overline{S}(L_{7})$.
	\end{proof}

	The observation made at the begining  of this section together with  Remark \ref{rk:join SF},  allows us to  characterize the supremum of two S-filters and the S-filter generated by a nonempty subset.\par 
	Let $F, G\in SF(L)$ and $\emptyset\neq X\subseteq L$.  Denote by  $F\underline{\vee} G$ the join of $F$ and $G$ in $SF(L)$ and  by  $S[X)$  the S-filter of $L$ generated by $X$. In the particular case  $X=\{a\}$,we write   $S[a)$ for the principal S-filter generated by $a$.

	\begin{proposition}\label{prop:Sprincp}  The following statements hold:
		\begin{enumerate}
			\item $S[ X)=\{x\in L\mid \exists x_{1},\ldots, x_{n}\in X, n\geq 1, x_{1}\overline{\sqcap}\ldots\overline{\sqcap} x_{n}\leq x\}$; in particular, 
			$$F\underline{\vee} G=\{x\in L\mid \exists u\in F, v\in G, u\overline{\sqcap} v\leq x\}.$$
			\item  $S[ a) =[a^{\Delta\Delta})$,  $S[ 0)=L$ and $S[ 1)=\{1\}$.
			\item  $S[ a) \underline{\vee} S[ b)= S[a\overline{\sqcap} b)$; in articular,  $S[a)\underline{\vee}S[a^{\Delta})=L.$\item  
			$S[ a) \cap S[ b) =S[ a\underline{\sqcup} b)$, in particular, $S[a)\cap S[a^{\Delta})=\{1\}$.
			
			\item 
			$a\leq b$ implies $S[ b)\subseteq S[ a)$.
		\end{enumerate}
	\end{proposition}	 
	\begin{proof}
		(1) Let $M:=\{x\in L\mid \exists x_{1}, \ldots, x_{n}\in X, n\geq 1, x_{1}\overline{\sqcap}\ldots\overline{\sqcap} x_{n}\leq x\}.$ Clearly  $X\subseteq M$, since 
		$x^{\Delta \Delta}\leq x$ for any   $x\in X$. Now,  let $x, y\in M$.  
		By remark \ref{lem:comp} (6), we have  $x\overline{\sqcap} y\in M$. 	
		If $x\leq y$, then by the transitivity of  $\leq$ and the definition of $M$,  we obtain $y\in M$. Therefore $M$ is an $S$-filter. 
		Now  let $F$ be an S-filter of $L$ containing $X$.  Let $x\in M$. Then there exist $x_{1}, \ldots, x_{n}\in X$ such that $$x_{1}\overline{\sqcap}\ldots\overline{\sqcap}x_{n}\leq x.$$
		Since  each $x_{i}\in X\subseteq F$,  and $F$ is closed under $\overline{\sqcap}$,  we deduce that $x\in F$. Thus, $M=S[X\rangle$.
		\par
		(2) Clearly,  $[a^{\Delta\Delta})\subseteq S[a)$. Let $z\in S[a)$. Taking $X=\{a\}$ in (1), one has $z\geq a^{\Delta\Delta}$, i.e $z\in[a^{\Delta\Delta})$. Therefore,  $S[a)\subseteq [a^{\Delta\Delta})$. Hence $[a^{\Delta\Delta})=S[a)$.\\
		(3) By (1) we have  $a\overline{\sqcap} b\in S[a)\underline{\vee} S[b)$, so $S[a\overline{\sqcap}b)\subseteq S[a)\underline{\vee} S[b)$. Let $z\in S[a\overline{\sqcap} b)$, then by (2), $z\geq (a\wedge b)^{\Delta\Delta}=a\overline{\sqcap} b)$. By (1), we conclude that $z\in S[a\overline{\sqcap}b)$. Thus $S[a\overline{\sqcap} b)=S[a)\underline{\vee} S[b)$. Taking $b=a^{\Delta}$,  we obtain the proof of the second part.\\
		(4) By  Remark \ref{lem:comp} (1)(ii),  we have $S[a\underline{\sqcup} b)\subseteq S[a)\cap S[b)$.  Let $z\in S[a)\cap S[b)$. Then, by (2), $a^{\Delta\Delta}\vee b^{\Delta\Delta}=a\underline{\sqcup}b\leq z$, so $z\in S[a\underline{\sqcup} b)$. Therefore $S[a)\cap S[b)\subseteq S[a\underline{\sqcup} b)$. Thus $S[a\underline{\sqcup})=S[a)\cap S[b)$.
		
		The second part follows from the fact that   $a\underline{\sqcup}a^{\Delta}=1$.\\
		(5) The last statement is traightforward.
	\end{proof}	 
	
	Let  $G\in F(\overline{S}(L))$.
	\begin{proposition}\label{prop: subext filter} The following statements hold.
		\begin{enumerate}\item  $G^{\star}=(F_{G})^{\star}=F_{G^{\overline{\star}}}$.
			\item   $G^{\overline{\star}~ \overline{\star}}=G$ if and only if $(F_{G})^{\star\star}=F_{G}$.
		\end{enumerate}
		
	\end{proposition}
	\begin{proof}

		(1) $G\subseteq F_{G}$, implies   $(F_{G})^{\star}\subseteq G^{\star}$.  Let $z\in G^{\star}$. Take  $u\in F_{G}$. Then  $u^{\Delta\Delta}\in G$,  and since $z\in G^{\star}$, it follows that $u^{\Delta\Delta\Delta}=u^{\Delta}\leq z$. Hence $z\in (F_{G})^{\star}$. Therefore $G^{\star}=(F_{G})^{\star}$.  \\
		$(\subseteq)$ Let $z\in F_{G^{\overline{\star}}}$. Then  $z^{\Delta\Delta}\in G^{\overline{\star}}$.  Let $x\in F_{G}$. Then $x^{\Delta\Delta}\in G$, and   $x^{\Delta}\leq z^{\Delta\Delta }\leq z$. Thus  $z\in (F_{G})^{\star}$ and  $F_{G^{\overline{\star}}}\subseteq (F_{G})^{\star}$.\\
		($\supseteq$) Conversely, let $z\in (F_{G})^{\star}$. For any $u\in G$, we also have  $u\in F_{G}$,  hence
		$u^{\Delta}\leq z^{\Delta\Delta}$. This means   $z^{\Delta\Delta}\in G^{\overline{\star}}$, so $z\in F_{G^{\overline{\star}}}$. Thus $F_{G^{\overline{\star}}}=(F_{G})^{\star}$.\\
		(2) Follows from (1) and Ramark \ref{prop: extfilter} (1). 
	\end{proof}
	
	The following theorem describes  how filters  of $L$ and filters of $\overline{S}(L)$ are related. Recall that $S(F(\overline{S}(L)))=\{G\in F(\overline{S}(L))\mid G^{\overline{\star}A \overline{\star}}=G\}$ and $S(SF(L))=\{F\in SF(L)\mid F^{\star\star}=F\}$.
	
	\begin{theorem}\label{theo: extiso} The following statements hold:
		\begin{enumerate}\item The map $\phi: F(\overline{S}(L))\to F(L), G\mapsto F_{G}$,  is  injective  and order-preserving. 	
			\item The algebras $(F(\overline{S}(\textbf{L}); \cap, \underline{\vee}, ^{\overline{\star}}, \{1\}, \overline{S}(L))$ and   $(SF(L); \cap, \underline{\vee}, ^{\star}, \{1\}, L)$ are isomorphic. 
			
			\item The algebras $(S(F(\overline{S}(L)); \cap, \underline{\vee}, ^{\overline{\star}}, \{1\}, \overline{S}(L))$ and $(S(SF(L); \cap, \underline{\vee}, ^{\star}, \{1\}, L)$ are isomorphic  ortholattices.
		\end{enumerate}
	\end{theorem}
	\begin{proof}
		(1) By   Proposition \ref{prop: extfilter} (1),  $\phi$ is well defined, injective and order preserving \par  
		(2) We have  $SF(L)=Im(\phi)$. Moreover,    $\phi$ is injective, order preserving, and satisfies $\phi(G^{\overline{\star}})=(F_{G})^{\star}$ (by   Proposition \ref{prop: subext filter}).\par
		(3) By Proposition \ref{prop: subext filter},  $G$ is $^{\overline{\star}~\overline{\star}}$ closed in $F(\overline{S}(L))$ if and only if $F_{G}$ is $^{\star\star}$ closed in $F(L)$. Since $\phi$ is injective, we  deduce that the corresponding  algebras  are isomorphic.  Furthermore,  since $S(SF(L))\subseteq SF(L)$ and $\phi: F(\overline{S}(L))\to SF(L)$ is surjective, we obtain the desired isomorphism. In addition $\phi(G^{\overline{\star}})=(\phi(G))^{\star}.$
	\end{proof}
	
	It is known that if $\textbf{L}$ is a complete lattice, then its set of principal filters is a sublattice of $(F(L); \subseteq))$ isomorphic to $\textbf{L}$( \cite{HP2002}). Similar results hold for ortholattices and WCL. For $a\in L$, we set $S[ a)^{\bot}=S[a^{\Delta})$ and consider the algebras $S\mathcal{F}_{p}(\textbf{L})=(SF_{p}(L); \underline{\vee}, \cap, ^{\bot}, S[ 1),S[ 0))$ of type (2,2,1,0,0),  where $SF_{p}(L)$ is the set of principal S-Filters of $L$.
	\begin{theorem}
		The algebra $S\mathcal{F}_{p}(\textbf{L})$ is an   ortholattice isomorphic to $\overline{S}(L)$.
	\end{theorem}
	\begin{proof}
		(1) By using (3), (4) and (5)  of Proposition \ref{prop:Sprincp} we have  $S[ a)^{\bot\bot}=S[ a\rangle$
		and the identities  $x\wedge x^{\bot}=0$ and $x\vee x^{\bot}=1$ hold in $SF_{p}(L)$. Hence $SF_{p}(L)$ is a bounded sublattice of $SF(L)$, and  therefore an ortholattice.
		Define  $$\Phi: \overline{S}(L)\to SF_{p}(L), a\mapsto S[ a).$$ 
		Clearly, $\Phi$ is well defined and compatible with operations, so it is a homomorphism of ortholattices.  By definition,  $\Phi$ is surjective. Now let  $a, b\in S(L)$ such that  $S[ a)=S[ b)$. We show that $a=b$. Since  $S[a)=S[b)$, we have $a\in S[ b)$, i.e   $b^{\Delta\Delta}=b\leq a$. Similarly,  $b\in S[ a)$ implies $a^{\Delta\Delta}=a\leq b$. Hence $a=b$. Thus $\Phi$ is injective. Consequently, $\Phi$ is an isomorphism of ortholattices.
	\end{proof}

	\section{Prime , maximal and primary filters in WDL}
	
	In Boolean algebras, distributive lattices, and double Boolean algebras, there are notions of ultrafilters, prime filters, and primary filters, which are closely related to the concepts of atoms, irreducible elements, prime elements, and primary elements. We study these notions  in WDL and investigate the relationship between the corresponding notions in a WDL $L$ and in its dual skeleton.
	\begin{proposition}Let $F$ be a filter of $L$.
		If $F$ is a prime  filter, then $F$ is a primary filter of $L$. Moreover, if $F$ is  a prime S-filter, then $F$ is maximal in $SF(L)$.
		
	\end{proposition}
	
	\begin{proof}
		Assume that $F$ is a prime filter  of $L$. For any $x\in L$, we have $x\vee x^{\Delta}=1\in F$. Since  $F$ is prime, it follows that $x\in F$ or $x^{\Delta}\in F$. Hence, $F$  is a primary filter. 
		Now suppose that $F$ is a  prime  S-filter. Then  $F$ is also a primary $S$-filter.  Let $x\in L\setminus F$. Then $x^{\Delta}\in F$. Moreover,  $$x\overline{\sqcap} x^{\Delta}=0\in S[ F\cup\{x\}\rangle_{F},$$ which shows that $F$ is maximal in $SF(L)$. 
	\end{proof}

	The correspondence betwen maximal filters of $\overline{S}(L)$ and maximal S-filters is established  in the following proposition. 
	
	\begin{proposition} The following statements hold:
		\begin{enumerate}
			\item $E$ is a maximal filter of $\overline{S}(L)$ if and only if  $F_{E}$ is a  maximal S-filter of $\textbf{L}$. 
			\item Each proper S-filter of $L$ is contained in a primary S-filter  of $L$. 	
		\end{enumerate}
	\end{proposition}
	\begin{proof}
		(1) Assume that  $E$ is a maximal filter  of $\overline{S}(L)$. Clearly,  $F_{E}$ is an S-filter of $L$ containing $E$. Let $H$ be an S-filter of $L$  such that $F_{E}\subseteq H$. Suppose $F_{E}\subsetneq H$ and let $x\in H\setminus F_{E}$.  Then $x^{\Delta\Delta}\not\in E$. Since $E$ is maximal, we have  $\langle E\cup\{x^{\Delta\Delta}\}\rangle=\overline{S}(L)$, so there exists $u\in E$ such that $u\overline{\sqcap} x^{\Delta\Delta}=0\in H$, which implies  $H=L$. Hence $F_{E}$ is a maximal S-filter of $L$.\\
		Conversely, assume  that $F_{E}$ is a primary S-filter of $\textbf{L}$, with $E\in F(\overline{S}(L))$.   Suppose   $E\nsubseteq H$ for some filter $H$ of $\overline{S}(\textbf{L})$. Let $x\in H\setminus E$. Then $x\not\in F_{E}$, so $x^{\Delta}\in E\subseteq H$. Hence $x\overline{\sqcap}x^{\Delta}=0\in H$, which implies   $H=\overline{S}(\textbf{L})$. Therefore, $E$ is a maximal filter of $\overline{S}(L)$.\par 		
		(2) Let $G$ be a proper $S$-filter of $L$. Since $G$ is proper, $0\not\in G$, and $E=G\cap \overline{S}(\textbf{L})$ is a proper filter of $\overline{S}(\textbf{L})$. Consider  $$\Psi:=\{H\subseteq \overline{S}(\textbf{L})\mid H~\text{is a proper filter of}~ \overline{S}(\textbf{L}), E\subseteq H\}.$$
		Assume that $\{H_{i}, i\in \Lambda\}\subseteq \Psi$ is totally ordered. Then $\underset{i\in\Lambda}{\bigcup}H_{i}$ is a proper filter of $\overline{S}(\textbf{L})$ containing $E$. Thus,  $\Psi$, ordered by inclusion, is inductive and by Zorn's Lemma there exists  a maximal filter $H$ of $\overline{S}(\textbf{L})$ containing $E$. Then  $F_{H}$ is a maximal S-filter containing $E$. By (1) $F_{H}$ is a primary and maximal S-filter.
	\end{proof}
	
	In the following theorem, we show  that primary filters and prime filters coincide in $\overline{S}(L)$, and that every maximal filter is a primary filter of $L$, although the converse does not hold.
	\begin{theorem}\label{theo: equivpp}
		Let $E$ be a filter of $\overline{S}(\textbf{L})$, $F$ a filter of $\textbf{L}$. The following statements hold.
		\begin{enumerate}
			\item $E$ is a primary filter of $\overline{S}(\textbf{L})$ if and only if $E$ is a prime filter of $\overline{S}(L)$.
			\item If $E$ is a primary filter of $\overline{S}(\textbf{L})$, then $E$ is a maximal filter of $\overline{S}(\textbf{L})$.
			
			\item If $F$ is a maximal filter of $\textbf{L}$, then $F$ is a primary filter of $\textbf{L}$. The converse is not true.
		\end{enumerate}
	\end{theorem}
	\begin{proof}
		(1) Assume that $E$ is a primary filter of $\overline{S}(L)$. We show that $E$ is a prime filter. Let $x, y\in \overline{S}(L)$ such that $x\vee y\in E$. Suppose by contradiction, that  $x\not\in E$ and $y\not\in E$. Since $E$ is primary,  $x^{\Delta}, y^{\Delta}\in E$. Then, using properties (4) and (5) of Proposition  \ref{prop: axiom}, we have 
		$$x^{\Delta}\overline{\sqcap} y^{\Delta}=(x^{\Delta}\wedge y^{\Delta})^{\Delta\Delta}=(x \vee y)^{\Delta}\in E.$$ Therefore,	$(x\vee y)\overline{\sqcap}((x\vee y)^{\Delta})=0\in E,$ which  contradicts the fact that  $E$ is a proper filter. Hence  $x\in E$ or $y\in E$, so  $E$ is prime.\\
		Conversely, assume that $E$ is a prime filter of $\overline{S}(L)$. For any $x\in \overline{S}(L)$, we have  $x\vee x^{\Delta}=1\in E$. Since   $E$ is prime, $x\in E$ or $x^{\Delta}\in E$. Thus, $E$ is primary. \par 
		(2) Assume  $E$ is a primary filter of $\overline{S}(L)$. Let $x\in \overline{S}(L)\setminus E$. Since $E$ is primary, $x^{\Delta}\in E$, and hence  
		$x\overline{\sqcap} x^{\Delta}=0\in [ E\cup\{x\})$. Therefore, $[E\cup\{x\})=\overline{S}(L)$, so  $E$ is  maximal. 
		\par
		\medskip
		(3) Assume that $F$ is a maximal filter of $L$. Clearly, $F$ is proper. Suppose, by contradiction, that  $F$ is not  primary. Then there exists $x\in L$ such that $x, x^{\Delta}\not\in F$. Since   $F$ is maximal,  we have $$[ F\cup\{x\})=[ F\cup\{x^{\Delta}\}) =L,$$
		so there exist $u, v\in F$ such that $u\wedge x=0$ and $v\wedge x^{\Delta}=0$. Then  $$(u\wedge v)\wedge x=(u\wedge v)\wedge x^{\Delta}=0.$$
		Using axiom (3), we get $((u\wedge v)\wedge x)\vee ((u\wedge v)\wedge x^{\Delta})=u\wedge v=0\in F,$ contradiction. Hence  $x\in F$ or $x^{\Delta}\in F$, so  $F$ is  primary.
		However,  consider the WDL $L_{6}$, $F=\{b,1\}$ is  primary but   not maximal.\end{proof}
	
	We denote by  
	$S\mathcal{F}_{Pr}(L)$  the set of all primary $S$-filter of $\textbf{L}$.
	Let $S\mathcal{F}(L)_{p}$ be the set of $S$- filters $F$ of $L$ such that $F\cap \overline{S}$ is a prime filter in $\overline{S}(L)$.
	\begin{proposition}\label{prop:iso-primary}
		We have 
		$
		S\mathcal{F}_{Pr}(L) = S\mathcal{F}(L)_{p}.$	
	\end{proposition}
	\begin{proof}
		Let $F$ be an $S$-primary filter of $L$. By (2) of Proposition \ref{prop: S-filt},  $F\cap \overline{S}(L)$ is a primary filter of $\overline{S}(L)$. Since  by Theorem \ref{theo: equivpp} (1), every  primary filter of $\overline{S}(L)$ is a prime filter, it follows that  $F\in S\mathcal{F}(L)_{p}$.
		For the reverse inclusion, let $F$ be an $S$-filter of $L$ such that $F\cap \overline{S}(L)=E$ is a prime filter of $\overline{S}(L)$. Then, by applying (2) of Proposition \ref{prop: S-filt}, we have  $F_{E}=F$. Hence,  $F$ is a primary $S$-filter of $L$ and  the equality holds. 
	\end{proof}
	The proposition \ref{prop:iso-primary}  implies that there is a one-to-one and onto correspondence between primary $S$-filters of $\mathbf{L}$ and prime filters of $\overline{S}(L)$.  
	
	An element $a \in L$ is called an \textbf{atom} (respectively, a \textbf{co-atom}) of $L$ if it covers $0$ (respectively, if it is covered by $1$). 
	\begin{remark}
		Let $L$ be a WDL, $a \in L$, and $F$ a filter of $L$.  
		\begin{enumerate}
			\item  $[a)$ of $L$ is maximal if and only if $a$ is an atom of $L$.
			
			\item  $S[a)$ is a maximal $S$-filter if and only if $a$ is an atom of $\overline{S}(L)$.
		\end{enumerate}
	\end{remark}

	\section{Congruences induced by S-filters in distributive WCL }

	The characterization of congruences provides tools that can simplify the study of WCLs. The description of subdirectly irreducible WDLs and simple WCLs depends on congruences generated by filters,  which offers natural way to  obtain the expected  characterizations.  $L$ denotes a WCL $(L; \vee, \wedge, ^{\Delta}, 0,1)$.
	An equivalence  relation  $\theta\subseteq  L^{2}$, is a congruence relation on  $L$ if $\theta$  is  compatible with $\vee,\wedge$ and  $^{\Delta}$. A congruence relation $\theta$ on $L$ is called \textbf{proper} if $\theta\neq L^{2}$. The set of all congruence relations on $L$ is denoted by $Con(L)$. For  $\theta\in Con(L)$ and $a\in L$, we write $$[a]_{\theta}:=\{x\in L\mid (a, x)\in\theta\}$$ for the equivalence class of $a$, in particular, $[1]_{\theta}$ is called, the \textbf{cokernel} of $\theta$.
	As stated in \cite{gratz2}, an equivalence relation $\theta$ on a lattice $(L; \vee, \wedge)$ is a congruence on $L$ if,  for any $x, y, z\in L$, $$(x, y)\in\theta~\text{ implies}~ (x\wedge z, y\wedge z), (x\vee z, y\vee z)\in\theta \quad (\ast).$$\par 
	Next, we consider  the binary relation   $\Phi$ defined on $L$ by 
	$$\Phi:=\{(a, b)\in L^{2}\mid   a^{\Delta}=b^{\Delta}\}.$$  It is easy to show that $\Phi\in Con(L)$. We call it \textbf{the determination congruence} of $L$.\\
	For a subset $S$ of $L$, $\theta, \psi\in Con(L)$,  $\theta_{S}$ denotes the restriction of $\theta$ to $S$ and  $\theta\circ \psi$ is the relational product of $\theta$ and $\psi$. An algebra $A$ is \textbf{congruence permutable} if  $\theta\circ \psi=\psi\circ\theta$ for all congruences $\theta$ and $\psi$ on $A$. A congruence relation is \textbf{regular} if $\theta=\psi$ whenever $\theta$ and $\psi$ have a class in common. As we shall see, the congruence  $\Phi$ plays a fundamental role in the study of distributive WCLs.
	It is known that $Con(L)$ forms a complete distributive lattice (\cite{gratz2}),  with a smallest element $$\Delta:=\{(a, a)\mid a\in L\}.$$
	In bounded distributive lattice, every filter is the cokernel of at least one congruence lattice.  This property does not carry over to WCLs. For instance, consider the WCL $L_{6}$ and the filter $F=\{a, b,1,v\}$, assume that $F=[1]_{\theta}$ for some congruence $\theta$ in $L_{6}$.
	As $a\in F$, we have  $(a, 1)\in\theta$, which implies  $(a^{\Delta\Delta}, 1)=(u, 1)\in\theta$ and therefore,  $u\in F$, contradiction. This observation leads us to characterize those filters that correspond to cokernels of congruences, and  describe the congruences generated by these filters.
	\begin{remark}
		Let $\theta\in Con(L)$. The  following statements hold.
		\begin{enumerate}
			\item (i) $F=[1]_{\theta}\in SF(L)$. (ii) If  $G=[1]_{\theta}\cap \overline{S}(L)$,  then $[1]_{\theta}=S[G)$.   	\item (i) If $(0, 1)\in\theta$, then $\theta=L^{2}$. (ii) 
			If  $\overline{S}(L)=\{0, 1\}$, then  $\theta=L^{2}$ or ( $[0]_{\theta}=\{0\}$ and $[1]_{\theta}=\{1\}$).	
		\end{enumerate}
	\end{remark}
	\begin{proof}	
		(1) For (ii), let $G=[1]_{\theta}\cap \overline{S}(L)$. Clearly,    $S[G)\subseteq [1]_{\theta}$. For the reverse inclusion,  let $z\in [1]_{\theta}$.  Then $(z, 1)\in \theta$, which  implies  $(z^{\Delta\Delta}, 1)\in\theta$, so  $z^{\Delta\Delta}\in G$. Since   $z^{\Delta\Delta}\leq z$, it follows from  (2) of Proposition \ref{prop: engendrer}  that $z\in S[G)$. Therefore the desired equality holds.\\
		(2)(i) Assume that $(0, 1)\in\theta$. Then  $0\in[1]_{\theta}$, $[1]_{\theta}\in F(L)$. Consequently,  $[1]_{\theta}=L$ and $\theta=L^{2}$. \\
		(ii)
		If $\overline{S}(L)=\{0, 1\}$, then the restriction of $\theta$ on its skeleton is a congruence represented by one of the following partitions: $\{\{0\},\{1\}\}, \{\{0,1\}\}$.
		\begin{itemize}
			\item For $(\{\{0\}, \{1\}\}$, we have 
			$[0]_{\theta}=I_{\{0\}}=\{0\}$ and $[1]_{\theta}=\{1\}$.
			\item In the second case  $(0, 1)\in \theta$, we have $\theta=L^{2}$.
		\end{itemize}
	\end{proof} 	
	\textbf{Problem}: given an S-filter $F$, how can we characterize the congruence generated by $F$?\\ We provide an answer to this question in the case of distributive WCLs.\\	
	From now, assume that $L$ is a distributive WCL. Let $F\in SF(L)$ and consider the binary relation
	$$\theta_{F}:=\{(x, y)\in L^{2}\mid \exists u\in F, x\vee u^{\Delta}=y\vee u^{\Delta}\}.$$ 
	\begin{theorem}\label{theo: cong induce}
		The binary relation $\theta_{F}$ is the least congruence relation on $L$ that collapse $F$. 
	\end{theorem}
	\begin{proof}	
		(1) 
		Clearly, $\theta_{F}$ is reflexive and symmetric. For  transitivity, let $x, y, z\in L$ such that $(x, y), (y, z)\in\theta_{F}$, then there exist $u, v\in F$ such that: 
		(i)~ $x\vee u^{\Delta}=y\vee u^{\Delta}$ and ~~(ii) $y\vee v^{\Delta}=z\vee v^{\Delta}$.\\ Since   $u\wedge v\in F$, we obtain  $$x\vee (u\wedge v)^{\Delta}=x\vee u^{\Delta}\vee v^{\Delta}=y\vee u^{\Delta}\vee v^{\Delta}=z\vee v^{\Delta}\vee u^{\Delta}=z\vee (u\wedge v)^{\Delta}.$$ Hence,  $(x, z)\in\theta_{F}$ and thus  $\theta_{F}$ is transitive.\\
		(2) Now we prove that $\theta_{F}$ is compatible with $\vee$, $\wedge$ and  $^{\Delta}$.\par 
		(2i)  Compatibility with $\wedge$ and $\vee$. Let $z\in L$ and assume  $(x, y)\in\theta_{F}$. Then there exists $u\in F$ such that $x\vee u^{\Delta}=y\vee u^{\Delta}$ ($\star_{1}$).
		For $\wedge$, using transitivity we have  $$(x\wedge z)\vee u^{\Delta}=(x\vee u^{\Delta})\wedge (z\vee u^{\Delta})=(y\vee u^{\Delta})\wedge (z\vee u^{\Delta})=(x\wedge z)\vee u^{\Delta}.$$ Therefore, $\theta_{F}$ is compatible with $\wedge$.
		For $\vee$,	by associativity of $\vee $ and by  ($\star_{1}$), we have  $$(x\vee z)\vee u^{\Delta}=(x\vee u^{\Delta})\vee z=(y\vee u^{\Delta})\vee z=(y\vee z)\vee u^{\Delta}.$$ It follows that $\theta_{F}$ is compatible with $\vee$.\par 
		(2ii) From the assumption, we have  $x\vee u^{\Delta}=y\vee u^{\Delta}$. Taking the meet with $u^{\Delta\Delta}$ on both  side of the last equality  and applying distributivity yields $(x\vee u^{\Delta})\wedge u^{\Delta\Delta}=(y\vee u^{\Delta})\wedge u^{\Delta\Delta}$ i.e $$(x\wedge u^{\Delta\Delta})\vee (u^{\Delta}\wedge u^{\Delta\Delta})=(y\wedge u^{\Delta\Delta})\vee (u^{\Delta\Delta}\wedge u^{\Delta\Delta}).$$ Since by (7) of Proposition \ref{prop: axiom} $u^{\Delta}\wedge u^{\Delta\Delta}=0$, this reduces to   $x\wedge u^{\Delta\Delta}=y\wedge u^{\Delta\Delta}.$ 
		Applying $^{\Delta}$ to both  sides gives  $x^{\Delta}\vee u^{\Delta\Delta\Delta}=y^{\Delta\Delta\Delta}\vee u^{\Delta\Delta\Delta}.$ As $u\in F$ and $u^{\Delta\Delta\Delta}=u^{\Delta}$,  by definition $(x^{\Delta}, y^{\Delta})\in\theta_{F}$. Hence $\theta_{F}$ is compatible with $^{\Delta}$.\par  From (1) and (2), it follows that $\theta$ is a congruence relation of $L$.\\
		Clearly, $F\subseteq [1]_{\theta_{F}}$. Let $x\in[1]_{\theta_{F}}$.  Then $(x, 1)\in\theta_{F}$, so there exists $u\in F$ such that $x\vee u^{\Delta}=1$. Taking the meet with  $u^{\Delta\Delta}$ on both  sides  and using distributivity, we   obtain  $x\wedge u^{\Delta\Delta}=u^{\Delta\Delta}$,   that is $u^{\Delta\Delta}\leq x\in F$. Since $u\in F$ and $F$ is an S-filter, it follows that $x\in F$. Hence  $[1]_{\theta_{F}}=F$.  Finally,  suppose that  $\psi$ is a congruence on $L$ such that $F\subseteq [1]_{\psi}$. 
		Let $(x, y)\in\theta_{F}$.  Then,  by definition,  there exists $u\in F$ such that $x\vee u^{\Delta}=y\vee u^{\Delta}.$ Taking the  meet with $u^{\Delta\Delta}$ on both sides of this  equality yields $$u^{\Delta\Delta}\wedge (x\vee u^{\Delta})=(x\wedge u^{\Delta\Delta})\vee (u^{\Delta\Delta}\wedge u^{\Delta})=x\wedge u^{\Delta\Delta}=y\wedge u^{\Delta\Delta}\quad (\star_{2}).$$ Since $(u^{\Delta\Delta}, 1)\in\psi$, the compatibility of $\psi$ with $\wedge$,  combined with ($\star_{2}$) implies $$(x\wedge u^{\Delta\Delta}, x), (y\wedge u^{\Delta\Delta}, y)\in\psi.$$ By transitivity of $\psi$, we obtain  $(x, y)\in\psi$.  Therefore,  $\theta_{F}\subseteq \psi$, and  hence $\theta_{F}=\Theta(F)$.
	\end{proof}
	
	Our next result gives a simple, explicit expression for the join of two congruences when one of them is a congruence generated by a filter. 
	\begin{theorem}
		Let $F\in SF(L)$ and $\Psi\in Con(L)$. Then $\Theta_{F}\vee \Psi=\theta_{F}\circ\Psi\circ\theta_{F}$.
	\end{theorem}
	\begin{proof}
		Note that if $(x, y)\in\theta_{F}$,  then there exists $f\in F$ such that $x\wedge f=y\wedge f$ ($\ast_{0})$.\\
		Clearly, one has $\theta_{F}\circ\Psi\circ\theta_{F}\subseteq \Psi\vee \theta_{F}.$
		For the reverse inclusion,  let us  first show that $$\Psi\circ\theta_{F}\circ\Psi\subseteq \theta_{F}\circ\Psi\circ \theta_{F}\quad (\ddagger\ddagger).$$
		Recall that $\Psi\vee \theta_{F}=\bigcup\{\alpha_{i_{1}}\circ\ldots\circ\alpha_{i_{k}}, k\geq 1, \alpha_{i_{l}}\in\{\Psi, \theta_{F}\}, \alpha_{i_{1}}=\Psi \} (\text{see \cite{Burris})}.$\\
		Note that,  using distributivity, $f\in F$,  one has  $(a\wedge f)\theta_{F} a$ for any $a\in L$ ($\star_{1})$.\\
		Let $(a,b)\in\Psi\circ\theta_{F}\circ\Psi$,  then there exist $u, v\in L$ such that $a\Psi u\theta_{F} v\Psi b.$ Since  $u\theta_{F}v$, by  ($\ast_{0})$ there exists $t\in F$ such that $u\wedge t=v\wedge t$~ ($\ast$). By   ($\ast_{1})$ we then have   $$a\theta_{F} a\wedge t,~ a\wedge t\Psi u\wedge t=v\wedge t,  v\wedge t\Psi b\wedge t, b\wedge t\theta_{F} b.$$ 
		It follows that  $$(a, a\wedge t)\in\theta_{F}, (a\wedge t, b\wedge t)\in\Psi,  
		(b\wedge t, b)\in\theta_{F}.$$ Therefore,  $(a, b)\in \theta_{F}\circ\Psi\circ\theta_{F}$, and ($\ddagger\ddagger$) is proved.\par  
		Clearly  $\Psi, \theta_{F}, \Psi\circ\theta_{F}, \theta_{F}\circ\Psi\subseteq \theta_{F}\circ\Psi\circ \theta_{F}.$
		Now we show by induction on $k\geq 1$ that $$\Psi\circ\theta_{F}\circ\Psi\circ\alpha_{1}\circ ..\circ \alpha_{k}\subseteq\theta_{F}\circ\Psi\circ\theta_{F},  \text{for}~\alpha_{j}\in\{\Psi, \theta_{F}\},  1\leq j\leq k.$$
		For $k=1$, if $\alpha_{1}=\Psi$, then by transitivity of $\Psi$ we are done. If $\alpha_{1}=\theta_{F}$, since $\theta_{F}$ is transitive and the composition  $\circ$ is isotone, composing each side of ($\ddagger\ddagger)$ on the right  by $\theta_{F}$ gives  $$\Psi\circ\theta_{F}\circ\Psi\circ\theta_{F}\subseteq \theta_{F}\circ\Psi\circ\theta_{F}.$$ Assume now  that the statement holds for some  $k\geq 1$.  
		A similar argument shows that  $$\Psi\circ\theta_{F}\circ\Psi\circ\alpha_{1}\circ\ldots\circ\alpha_{k}\circ\alpha_{k+1}\subseteq \theta_{F}\circ\Psi\circ\theta_{F}.$$ 
		Hence, for any $k\geq 1$, we have  $$\theta_{F}, \Psi, \theta_{F}\circ\Psi, \Psi\circ\theta_{F}, \Psi\circ\theta_{F}\circ\Psi\circ\alpha_{1}\circ\ldots\circ\alpha_{k}\circ\alpha_{k+1}\subseteq \theta_{F}\circ\Psi\circ\theta_{F}.$$  Thus $\Psi\vee\theta_{F}=\theta_{F}\circ\Psi\circ\theta_{F}.$
	\end{proof}

	\begin{theorem}
		\begin{enumerate}
			\item $\Phi$ is the largest congruence $\eta$  in $L$ with $[0]_{\eta}=\{0\}$.
			\item $L$ is regular if and only if $\Phi=\Delta$.
		\end{enumerate}
		
	\end{theorem}
	
	\begin{proof}
		(1)  Let $\eta$ be a congruence on $L$ such that  $[0]_{\eta}=\{0\}$.  Suppose  $a\cong b(\eta)$. By compatibility with $^{\Delta}$, we have   $(a^{\Delta}, b^{\Delta})\in (\eta)$, therefore, for any $x\in L$, ~$(x\wedge  a^{\Delta},  x\wedge b^{\Delta})\in\eta$  $\quad (\star).$  Taking $x=b^{\Delta\Delta}$ in ($\star$), we obtain  $(b^{\Delta\Delta}\wedge b^{\Delta}=0,  b^{\Delta\Delta}\wedge a^{\Delta})\in (\eta).$ Since  $[0]_{\eta}=\{0\}$, it follows that  $b^{\Delta\Delta}\wedge  a^{\Delta}=0$,  by (8) of Proposition \ref{prop: axiom},  we deduce that   $b^{\Delta}\geq a^{\Delta}$. Similarly, taking $x=a^{\Delta\Delta}$ in $(\star$), one get $b^{\Delta}\leq a^{\Delta}$.  Hence $a^{\Delta}=b^{\Delta}$, that is $(a,b)\in\Phi$.   Hence $\eta\subseteq\Phi$. 
		\par 
		(2) For the  direct implication, if $L$ is regular, it is clear that $\Delta=\Phi$.  Indeed,  $[0]_{\Phi}=\{0\}$ and $\Delta$  shares  a class with $\Phi$;  by regularity, it follows that $\Phi=\Delta$.\\ Conversely,  assume that $\Delta=\Phi$.  Let $\theta$  be a congruence on $L$ and  $\{b\}$ be a single  class of $\theta$. Let $I:=[0]_{\theta}$. For any   $a\in I$, we have $a\cong 0$, hence  $(a^{\Delta\Delta}, 0)\in \theta$. It follows that   $(b\vee a^{\Delta\Delta},  b)\in \theta$, which implies   $a^{\Delta\Delta}\vee b=b$, i.e  $ a^{\Delta\Delta}\leq b$. On the other hand  $(a^{\Delta}, 1)\in \theta$, implies $(a^{\Delta}\wedge b, b)\in \theta$, hence $a^{\Delta}\wedge b=b$, so $b\leq a^{\Delta}$.  Combining these inequality  we have  $a^{\Delta\Delta}\leq b\leq a^{\Delta}$. Since $a^{\Delta\Delta}\vee a^{\Delta}=1$, it follows that $a^{\Delta}=1=0^{\Delta}$, so   $(a, 0)\in\Phi=\Delta$, and $a=0$. Thus $I=[0]_{\theta}=\{0\}$. By part (1) $\Phi$ is the largest congruence $\eta$ in $L$ with $[0]_{\eta}=\{0\}$. Therefore, $\theta=\Delta$. Hence $L$ is regular. 
	\end{proof}
	
	\begin{theorem}
		The congruences generated by S-filters commute.
	\end{theorem}
	\begin{proof}
		Let $F_{1}, F_{2}$ be two S-filters of $L$. Assume that $x\theta_{F_{1}} z\theta_{F_{2}} y$, then there exist $f_{1}\in F_{1}$ and $f_{2}\in F_{2}$ such that $$x\vee f_{1}^{\Delta}=z\vee f_{1}^{\Delta} \quad (\dagger),~ \text{and}\quad z\vee f_{2}^{\Delta}=y\vee f_{2}^{\Delta}\quad  (\dagger_{2}).$$ Let $\omega=(x\vee f_{2}^{\Delta})\wedge (y\vee f_{1}^{\Delta})$. Note that $$(x\vee f_{2}^{\Delta}, x)\in ~ \theta_{F_{2}})~\quad   (\dagger_{3}), \text{and} \quad (y\vee f_{1}^{\Delta},  y)\in \theta_{F_{1}}~\quad  (\dagger_{4}).$$ 
		From  $y\cong z\theta_{F_{2}}$, it follows that   $(y\vee f_{1}^{\Delta}, z\vee f_{1}^{\Delta})\in \theta_{F_{2}}.$ Using ($\dagger_{4}$), we obtain  $$((y\vee f_{1}^{\Delta})\wedge (x\vee f_{2}^{\Delta}), (z\vee f_{1}^{\Delta})\wedge x=x)\in \theta_{F_{2}}.$$  Therefore, $(\omega,  x)\in\theta_{F_{2}}$.  Similarly, from $(x, z)\in \theta_{F_{1}}$, we get  $(x\vee f_{2}^{\Delta},  z\vee f_{2}^{\Delta})\in \theta_{F_{1}}.$ Using ($\dagger_{3}$), this gives $(\omega, (z\vee f_{2}^{\Delta})\wedge y=y)\in \theta_{F_{1}}.$ Thus  we have both  $(\omega, y)\in \theta_{F_{1}})$ and $(x, \omega)\in  \theta_{F_{2}}$. Hence $\theta_{F_{1}}\circ\theta_{F_{2}}\subseteq\theta_{F_{2}}\circ \theta_{F_{1}}.$ A symmetric  argument shows that $\theta_{F_{2}}\circ\theta_{F_{1}}\subseteq\theta_{F_{1}}\circ\theta_{F_{2}}.$\\ Therefore,  the desired equality holds.
	\end{proof}
	\begin{remark}\label{lem: inj}
		Let  $F, G\in SF(L)$. The following statement holds:
		
		$$F\subseteq G~ \text{ iff}\quad  \theta_{F}\subseteq  \theta_{G}.$$ 
	\end{remark}
	
	
	
	We now provide some examples of  congruences $\theta_{F}$ with the corresponding S-filters.

	\begin{example}
		Consider the  WCL $L_{6}$. Let $\theta$ (rep. $\Psi$) be the equivalence relation defined on $L_{6}$ by its classes $\{0,u\}, \{b,1\}, \{a,v\}$ (resp.  $\{0,b,v\}, \{u,a,1\}$). Then:
		\begin{enumerate}
			\item   
			 $\theta$ is a congruence relation on $L_{6}$ with  
			the corresponding S-filter  $F=[1]_{\theta}$; 
			
			\item  
			$\Psi$ is a congruence relation on $L_{6}$, and the corresponding S-filter is $F=\{u,a,1\}$.
		\end{enumerate}
	\end{example}

	\begin{proposition}(G. Birkhoff, 1944). An algebra $\mathcal{A}$ is subdirectly irreducible if and only if it is nontrivial and $\cap\{\theta\in Con(\mathcal{A})\mid \theta\neq\Delta_{A}\}\neq\Delta_{A}$, which means that the lattice $Con(\mathcal{A})$ has just one atom.	
	\end{proposition}
	
	\begin{lemma}  Let  $\theta, \beta\in Con(L)$. The following statements holds. 
		
		\begin{enumerate}\item $\Phi$ is the largest  congruence $\theta$  on $L$ verifying $\theta_{\overline{S}(L)}=\Delta_{\overline{S}(L)}.$
			\item If $\Phi=\Delta$, then  $[1]_{\theta}=\{1\}$ iff  $\theta=\Delta_{L}$.		
		\end{enumerate}
	\end{lemma}
	\begin{proof}
		
		(1) Let $\theta\in Con(L)$.  Let $(x, y)\in\theta$, then $(x^{\Delta\Delta}, y^{\Delta\Delta})\in \theta_{S(L)}=\Delta_{\overline{S}(L)}$, so $x^{\Delta\Delta}=y^{\Delta\Delta}$. Hence  $(x, y)\in\Phi$ and therefore, $\theta\subseteq \Phi$.\par 
		(2) Assume that $\Phi=\Delta$ and let $\theta\in Con(L)$. Suppose that $[1]_{\theta}=\{1\}$. Then one can show that $[0]_{\theta}=\{0\}$. Thus  $\theta\subseteq\Phi=\Delta$. Hence $\theta=\Delta$.
	\end{proof}
	\begin{proposition}\label{pro: irre}
		If  $\Delta=\Phi$, then  $\bigcap\{ \theta: \theta\in Con(L)\setminus\{\Delta\}\}=\bigcap\{\theta_{F}\mid F\in SF(L)\setminus\{1\}\}\}.$
	\end{proposition}
	
	\begin{proof}	
		Recall that $\Delta = \theta_{\{1\}}$ and $L^{2} = \theta_{L}$. Set        $$\gamma=\bigcap \{\theta: \theta \in Con(L)\setminus\{\Delta\}\},  ~\text{and}~ \psi=\bigcap\{ \theta_{F}: F\in  SF(L)\setminus \{\{1\}\}\}.$$   Now we  show that $\gamma=\psi$. Clearly $\gamma\subseteq\psi$. For the reverse inclusion, let $(x, y)\in\psi$. Then,  for any S-filter $F\neq \{1\}$ on $L$, we have 
		$(x, y)\in \theta_{F}~ (\star).$
		Let $\theta$ be a congruence relation on $L$ with  $\theta\neq\Delta$. Since $\Delta=\Phi$ and $\theta$ is proper, it follows that  $F=[1]_{\theta}\neq \{1\}$
		is a proper S-filter of $L$. Hence  $(x, y)\in\theta_{F}\subseteq \theta$, so $(x,y)\in\theta$ for every proper congruence $\theta$ on $L$. Therefore,  the desired   equality holds and we conclude that  $\gamma=\psi$. 
	\end{proof}	
	\begin{theorem}
		If $\Delta=\Phi$, then $Con(L)$ and $NF(L)$ are isomorphic lattices. 
	\end{theorem}
	\begin{proof}
		Clearly $NF(L)$ embeds in $Con(L)$ by Remark \ref{lem: inj}. Let $\Psi\in Con(L)$ with $\Psi\neq \Delta$. Then $F=[1]_{\Psi}\neq\{1\}$ by assumption on $L$.  Let $\theta=\theta_{F}$. One has $[1]_{\theta_{F}}=F=[1]_{\Psi}$, as $L$ is a regular WCL, we deduce that $\Psi=\theta_{F}$.  Therefore, $h$ is surjective. Hence $Con(L)$ and $SF(L)$ are isomorphic lattices.
	\end{proof}
	\begin{theorem} If 
		$\Phi=\Delta_{L}$, then the following statements hold.
		\begin{enumerate}\item  $L$ is subdirectly irreducible if and only if $SF(L)$ has a unique atom.\item  $L$ is simple if and only if $\mid SF(L))\mid =2$. \end{enumerate} 
	\end{theorem}
	\begin{proof}
		
		Assume that $\Delta=\Phi$. Then $Con(L)=\{\theta_{F}\mid F\in NF(L)\}.$ If there exists a least S- filter $F\neq \{1\}$ of $L$ such that $F\subseteq G$ for all $G\in NF(L), G\neq \{1\}$, then $\theta_{F}=\cap ( Con(L)\setminus\{\Delta\})$ is the monolith congruence in $L$ and $L$ is subdirectly irreducible. Then $F=N[e)$ for some $e\in L$.\\  Conversely, assume that $L$ is subdirectly irreducible, then there exists a least congruence $\mu$ on $L$ such that  $\mu\neq\Delta$ and   $$\mu=\cap\{\theta\mid \theta\in Con(L)\setminus\{\Delta\}\}\}=\cap\{\theta_{F}\mid F\in NF(L), F\neq \{1\}\}$$ (by Proposition \ref{pro: irre}), since $\Delta=\Phi$, $\mu=\theta_{[1]_{\mu}}$ and $[1]_{\mu}$ is the least    S-filter of $L$  distinct from $\{1\}$.		
	\end{proof}

	\section{Conclusion}
	In this work,  we have  provided a dual weak complementation operation in the set of filters of a given WDL derived from $^{\Delta}$. We have also proved that each filter of a skeleton of a WDL  is   a base of some filter in $L$ . The collection of such filters forms  a complete lattice isomorphic to the filter lattice of the dual skeleton of $L$. S-primary filters are characterized and it is proved that each S-primary filter is maximal in the lattice of S-filters and corresponds to a prime filter in the dual skeleton $\overline{S}(L)$.  S-filters introduced are used to construct induced congruences in distributive weakly complemented lattices. Simple and subdirectly irreducible  weakly complemented lattices are characterized in term of S-filters.  Finally, maximal  principal filters  are charaterized in term of atom  of $L$. Our next step to investigate WDL is concerned with the topological representation of WDL and  congruences generated by filters  in WDL.\\

\end{document}